# A note on some versions of ♣

Mirna Džamonja
Mathematics Department
University of Wisconsin-Madison
Madison, WI 53706, USA

dzamonja@math.wisc.edu

Saharon Shelah
Mathematics Department
Hebrew University of Jerusalem
91904 Givat Ram, Israel

shelah@sunset.huji.ac.il

September 12, 1997

**Abstract**

This note gives two results on guessing unbounded subsets of $\lambda^+$. The first is a positive result and applies to the situation of $\lambda$ regular, while the second is a negative consistency result which applies to the situation of $\lambda$ singular. Both results are connected to an earlier result of Džamonja-Shelah (see Fact 0.2 of this paper) in which they showed that a certain version of ♣ holds at a successor of singular just in $ZFC$. The first result here shows that the result of Fact 0.2 can to a certain extent be extended to the successor of regular. The negative result here gives limitations to the extent to which one can hope to extend Fact 0.2. [1]

---

[1]This note is numbered F227(10/97) in Saharon Shelah's list of publications. It is presently available only as an electronic publication.

# 0   Introduction and background

We consider possible improvements of a result from [DjSh 545] about the $ZFC$ existence of a certain version of ♣. The result in question is quoted as Fact 0.2 below, and can be summarized as having shown that at the succesor of singular a certain version of ♣ always hold. One direction that can be suggested, is to consider the successor of regular. A positive result in this direction is given in §1 here.

Another direction is to consider $ZFC$ improvements to Fact 0.2 at the successor of singular. In §2 we show that there are some non-obvious limitations, by giving a consistency result regarding successor of singular strong limit in which a negation of guessing is obtained.

We now recall the relevant result from [DjSh 545].

**Definition 0.1.** Suppose that $\lambda$ is a cardinal. $\clubsuit^*_{-\lambda}(\lambda^+)$ is the statement saying that there is a sequence $\langle \mathcal{P}_\delta : \delta \text{ limit } < \lambda^+ \rangle$ such that

(i) $\mathcal{P}_\delta$ is a family of $\leq |\delta|$ unbounded subsets of $\delta$,

(ii) For $a \in \mathcal{P}_\delta$ we have $\text{otp}(a) < \lambda$,

(iii) For all $X \in [\lambda^+]^{\lambda^+}$, there is a club $C$ of $\lambda^+$ such that for all $\delta \in C$ limit, there is $a \in \mathcal{P}_\delta$ such that
$$\sup(a \cap X) = \delta.$$

**Fact 0.2 (Džamonja-Shelah).** [DjSh 545] If $\aleph_0 < \kappa = \text{cf}(\lambda) < \lambda$, then $\clubsuit^*_\lambda(\lambda^+)$.

In §1 we show that Fact 0.2 can be to a certain extent extended to the successor of regular, by proving the following (see below for the relevant notation):

**Theorem 0.3.** 1.1

(1) Suppose that

(a) $\lambda = \mathrm{cf}(\lambda) > \theta_1^+ > \theta_1 \geq \theta = \mathrm{cf}(\theta) > \aleph_1$ and $\lambda \geq 2^{\aleph_0}$,

(b) $S^* \subseteq S_\theta^{\lambda^+}$ is stationary, moreover
$$S_1 \stackrel{\mathrm{def}}{=} \{\delta < \lambda^+ : \mathrm{cf}(\delta) = \theta_1^+ \ \& \ S^* \cap \delta \text{ is stationary}\}$$
is stationary.

Then there is a stationary $S' \subseteq S^*$ and $\langle E_\delta : \delta \in S' \rangle$ such that

(i) $E_\delta$ is a club of $\delta$ with $\mathrm{otp}(E_\delta) \leq \lambda^\omega \times \theta_1^+$,

(ii) for every unbounded $A \subseteq S_\lambda^{\lambda^+}$, for stationarily many $\delta \in S'$, we have
$$\delta = \sup(A \cap \mathrm{nacc}(E_\delta)).$$

(2) We can omit the assumption of $\lambda \geq 2^{\aleph_0}$ if $S^* = S_\theta^{\lambda^+}$.

Note that this result is in some sense complementary to the "club guessing" results of Shelah, because here we are guessing unbounded subsets of $\lambda^+$ which are not necessarily clubs, but on the other hand, there are limitations on the cofinalities.

In §2, modulo the existence of a supercompact (Theorem 2.1), it is shown here that it is consistent that there is $\lambda$ a strong limit singular of cofinality $\omega$, such that $2^\lambda > \lambda^+$ and the following negation of guessing holds:

There is a function $f : \lambda^+ \to \omega$ such that for every $\mathcal{P} \subseteq [\lambda^+]^\omega$ of cardinality $< 2^\lambda$, for some $X \in [\lambda^+]^{\lambda^+}$ we have

(i) $(\forall i < \omega)[|X \cap f^{-1}(\{i\})| = \lambda^+]$,

(ii) If $a \in \mathcal{P}$, then $\sup(\mathrm{Rang}(f \upharpoonright (a \cap X))) < \omega$.

We finish this introduction by recalling some notation and facts which will be used in the following sections.

**Notation 0.4.** (1) Suppose that $\kappa = \mathrm{cf}(\kappa) < \delta$. We let
$$S_\kappa^{\prime\delta} \stackrel{\mathrm{def}}{=} \{\alpha < \delta : \mathrm{cf}(\alpha) = \kappa\}.$$

(2) Suppose that $C \subseteq \alpha$. We let
$$\mathrm{acc}(C) \stackrel{\mathrm{def}}{=} \{\beta \in C : \beta = \sup(C \cap \beta)\},$$
and $\mathrm{nacc}(C) \stackrel{\mathrm{def}}{=} C \setminus \mathrm{acc}(C)$.

**Definition 0.5.** Suppose that $\lambda \geq \aleph_1$ and $\gamma$ is an ordinal, while $A \subseteq \lambda^+$. For $S \subseteq \lambda^+$, we say that *$S$ has a square of type $\leq \gamma$ nonaccumulating in $A$* iff there is a sequence $\langle e_\alpha : \alpha \in S \rangle$ such that

(i) $\beta \in e_\alpha \implies \beta \in S \ \& \ e_\beta = e_\alpha \cap \beta$,

(ii) $e_\alpha$ is a closed set,

(iii) If $\alpha \in S \setminus A$, then $\alpha = \sup(e_\alpha)$,

(iv) $\mathrm{otp}(e_\alpha) \leq \gamma$.

**Fact 0.6 (Shelah).** [[Sh 351]§4, [Sh -g]III§2] Suppose that
$$\lambda = \mathrm{cf}(\lambda) > \aleph_1, \kappa = \mathrm{cf}(\kappa).$$
Further suppose that $S \subseteq S^\lambda_\kappa$ is stationary. <u>Then</u> there is $S_1 \subseteq \lambda$ on which there is a square of type $\leq \kappa$, nonaccumulating on $A$=the successor ordinals, and $S_1 \cap S$ is stationary.

**Remark 0.7.** In the proof of Fact 0.6 we can replace $A$=the successor ordinals with $A = S^\lambda_\sigma$ for any $\sigma = \mathrm{cf}(\sigma) < \kappa$.

**Definition 0.8 (Shelah).** [Sh -g] Suppose that $\delta < \lambda$ and $e \subseteq \delta$, while $E \subseteq \lambda$. We define
$$\mathrm{gl}(e, E) \stackrel{\mathrm{def}}{=} \{\sup(\alpha \cap E) : \alpha \in e \ \& \ \alpha > \min(E)\}.$$

**Observation 0.9.** Suppose that $e$ and $E$ are as in Definition 0.8, and both $e$ and $E \cap \delta$ are clubs of $\delta$. <u>Then</u>, observe that $\mathrm{gl}(e, E)$ is a club of $\delta$ with $\mathrm{otp}(\mathrm{gl}(e, E)) \leq \mathrm{otp}(e)$.

If $e$ is just closed in $\delta$, and $E \cap \delta$ is a club of $\delta$ <u>then</u> $\mathrm{gl}(e, E)$ is closed and $\mathrm{otp}(\mathrm{gl}(e, E)) \leq \mathrm{otp}(e)$.

**Fact 0.10 (Shelah).** [Sh 355] Suppose that $\mathrm{cf}(\kappa) = \kappa < \kappa^+ < \mathrm{cf}(\lambda) = \lambda$. Further suppose that $S \subseteq S_\kappa^\lambda$ is stationary and $\langle e_\delta : \delta \in S \rangle$ is a sequence such that each $e_\delta$ is a club of $\delta$. <u>Then</u> there is a club $E^*$ of $\lambda$ such that the sequence
$$\bar{c} = \langle c_\delta \stackrel{\mathrm{def}}{=} \mathrm{gl}(e_\delta, E^*) : \delta \in S \cap E^* \rangle$$
has the property that for every club $E$ of $\lambda$, there are stationarily many $\delta$ such that $c_\delta \subseteq E$.

**Observation 0.11.** Suppose that $\mathrm{cf}(\kappa) = \kappa < \kappa^+ < \mathrm{cf}(\lambda) = \lambda$ and that $S_1 \subseteq S_\kappa^\lambda$ is stationary, while $A = S_\sigma^\lambda$ for some $\sigma = \mathrm{cf}(\sigma) < \kappa$, possibly $\sigma = 1$. <u>Then</u> there is stationary $S_2 \subseteq \lambda$ and a square $\langle e_\delta : \delta \in S_2 \rangle$ of type $\leq \kappa$ nonaccumulating in $A$, such that $S_1 \cap S_2$ is stationary and

$$E \text{ a club of } \lambda \Longrightarrow \{\delta \in S_1 \cap S_2 : e_\delta \subseteq E\} \text{ is stationary.}$$

[Why? By Fact 0.6, there is $S_3 \subseteq S_1$ with a square $\langle e_\alpha : \alpha \in S_3 \rangle$ of type $\leq \kappa$ nonaccumulating in $A$, and that $S_1 \cap S_3$ is stationary. By Fact 0.10, there is club $E^*$ of $\lambda$ as in the conclusion of Fact 0.10, with $S_3 \cap S_1$ in place of $S$. Now, letting

$$S_2 \stackrel{\mathrm{def}}{=} \{\sup(\alpha \cap E^*) : \alpha \in \bigcup_{\delta \in S_3} e_\delta \cup \{\delta\} \ \& \ \alpha > \min(E^*)\},$$

and for $\delta \in S_2$, letting $c_\delta \stackrel{\mathrm{def}}{=} \mathrm{gl}(e_\delta, E^*)$, we observe that $S_2 \cap S_1$ is stationary (as $S_2 \cap S_1 \supseteq S_1 \cap S_3 \cap \mathrm{acc}(E^*)$), and $\langle c_\delta : \delta \in S_2 \rangle$ is a square of type $\leq \kappa$ nonaccumulating in $A$, while

$$E \text{ a club of } \lambda \Longrightarrow \{\delta \in S_1 \cap S_2 : c_\delta \subseteq E\} \text{ is stationary.}]$$

# 1  A $ZFC$ version of ♣

**Theorem 1.1.**  (1) Suppose that

  (a) $\lambda = \mathrm{cf}(\lambda) > \theta_1^+ > \theta_1 \geq \theta = \mathrm{cf}(\theta) > \aleph_1$ and $\lambda \geq 2^{\aleph_0}$,

(b) $S^* \subseteq S_\theta^{\lambda^+}$ is stationary, moreover

$$S_1 \stackrel{\text{def}}{=} \{\delta < \lambda^+ : \text{cf}(\delta) = \theta_1^+ \ \& \ S^* \cap \delta \text{ is stationary}\}$$

is stationary.

Then there is a stationary $S' \subseteq S^*$ and $\langle E_\delta : \delta \in S' \rangle$ such that

(i) $E_\delta$ is a club of $\delta$ with $\text{otp}(E_\delta) \leq \lambda^\omega \times \theta_1^+$,

(ii) for every unbounded $A \subseteq S_\lambda^{\lambda^+}$, for stationarily many $\delta \in S'$, we have
$$\delta = \sup(A \cap \text{nacc}(E_\delta)).$$

(2) We can omit the assumption of $\lambda \geq 2^{\aleph_0}$ if $S^* = S_\theta^{\lambda^+}$.

**Proof.** (1) Let $S_0 \stackrel{\text{def}}{=} S_\theta^{\lambda^+}$ and let $A^* \stackrel{\text{def}}{=} S_{\aleph_1}^{\lambda^+}$.

By Observation 0.11, there is a $S_2 \subseteq S_{\leq \theta_1^+}^{\lambda^+}$ such that there is a square $\langle e_\delta : \delta \in S_2 \rangle$ of type $\leq \theta_1^+$ nonaccumulating in $A^*$, the set $S_1 \cap S_2$ is stationary, and, moreover, for every $E$ a club of $\lambda^+$, the set $\{\delta \in S_1 \cap S_2 : e_\delta \subseteq E\}$ is stationary. [Why can we assume that $S_2 \subseteq S_{\leq \theta_1^+}^{\lambda^+}$? Just throw away the elements of higher cofinality.]

Let $S' \stackrel{\text{def}}{=} S^* \cap S_2$, so stationary.

**Claim 1.2.** There is a function $g : S' \to \omega$ such that for every club $E$ of $\lambda^+$, there are stationarily many $\delta \in S_1 \cap S_2$ such that $e_\delta \subseteq E$ and

$$(\forall n < \omega)[E \cap \delta \cap g^{-1}(\{n\}) \text{ is stationary in } \delta].$$

**Proof of the Claim.** For $\delta \in S^*$, we choose a sequence $\langle \xi_{\delta,i} : i < \theta \rangle$ increasing with limit $\delta$, and such that $\xi_{\delta,i} \in e_\delta$ and $\text{otp}(e_{\xi_{\delta,i}})$ depends only on $i$ and $\text{otp}(e_\delta)$, but not on $\delta$. For each $i < \theta$, we define a function $h_i : S' \to \theta_1^+$ by letting

$$h_i(\alpha) \stackrel{\text{def}}{=} \text{otp}(e_{\xi_{\alpha,i}}).$$

For $\beta < \theta_1^+$, let $h_i^\beta : S' \to \beta + 1$ be given by $h_i^\beta(\alpha) \stackrel{\text{def}}{=} \min\{h_i(\alpha), \beta\}$.

**Subclaim 1.3.** For some $\beta = \beta^*, i = i^*$, letting $h = h_i^\beta$ we have that for every club $E$ of $\lambda^+$, there are stationarily many $\delta \in S_1 \cap S_2$ such that $e_\delta \subseteq E$ and
$$|\{\zeta \leq \beta : E \cap \delta \cap h^{-1}(\{\zeta\}) \cap S' \text{ is stationary in } \delta\}| \geq \theta_1.$$

**Proof of the Subclaim.** Suppose the subclaim is not true, so for every $\beta < \theta_1^+$, and $i < \theta$ we can find clubs $E_{\beta,i}^0$ and $E_{\beta,i}^1$ of $\lambda^+$ such that $[\delta \in S_1 \cap S_2 \cap E_{\beta,i}^1 \ \& \ e_\delta \subseteq E_{\beta,i}^0] \implies$

$$|\{\zeta \leq \beta : E_{\beta,i}^0 \cap \delta \cap (h_i^\beta)^{-1}(\{\zeta\}) \cap S' \text{ is stationary in } \delta\}| < \theta_1.$$

Let $E^* \overset{\text{def}}{=} \text{acc}(\bigcap_{\beta,i} E_{\beta,i}^0 \cap \bigcap_{\beta,i} E_{\beta,i}^1)$, so a club of $\lambda^+$. Let $\delta \in E^* \cap S_1 \cap S_2$ be such that $e_\delta \subseteq E^*$ and $\delta > \theta_1^+ + 1$. Hence, for all $\beta < \theta_1^+$ and $i < \theta$,

$$|\{\zeta \leq \beta : E^* \cap \delta \cap (h_i^\beta)^{-1}(\{\zeta\}) \cap S' \text{ is stationary in } \delta\}| < \theta_1.$$

Fix $i < \theta$. Note:

$$\zeta < \beta_1 \leq \beta_2 \implies (h_i^{\beta_1})^{-1}(\{\zeta\}) = (h_i^{\beta_2})^{-1}(\{\zeta\}) = h_i^{-1}(\{\zeta\}).$$

By induction on $\varepsilon < \theta_1$ we try to choose $\beta_\varepsilon, \zeta_\varepsilon$ such that

(i) $\beta_{\varepsilon+1} > \zeta_\varepsilon > \beta_\varepsilon$,

(ii) $(h_i^{\beta_{\varepsilon+1}})^{-1}(\{\zeta_\varepsilon\}) \cap S' \cap E^*$ is stationary in $\delta$.

The induction must stop somewhere, as otherwise, taking

$$\beta^* \overset{\text{def}}{=} \sup\{\beta_\varepsilon : \varepsilon < \theta_1\},$$

we get $\beta^* < \theta_1^+$, yet

$$\{\zeta \leq \beta^* : (h_i^{\beta^*})^{-1}(\{\zeta\}) \cap S' \cap E^* \text{ is stationary in } \delta\} \supseteq \{\zeta_\varepsilon : \varepsilon < \theta_1\},$$

a contradiction. Hence, there is $\beta^i$ such that for no $\theta_1^+ > \beta > \zeta > \beta^i$ do we have
$$(h_i^\beta)^{-1}(\{\zeta\}) \cap S' \cap E^* \text{ is stationary in } \delta.$$

Let $\beta^* \stackrel{\text{def}}{=} \sup\{\beta^i : i < \theta\}$, so $\beta^* < \theta_1^+$ and for no $\theta_1^+ > \beta > \zeta > \beta^*$ and $i < \theta$ do we have that $(h_i^\beta)^{-1}(\{\zeta\}) \cap S' \cap E^*$ is stationary in $\delta$. Note that for every such $\zeta, \beta$ we have $(h_i^\beta)^{-1}(\{\zeta\}) = h_i^{-1}(\{\zeta\})$, for every $i < \theta$. Hence, for every $i < \theta$ there is $\zeta_i$ such that

$$\zeta \geq \zeta_i \Longrightarrow h_i^{-1}(\{\zeta\}) \cap S' \cap E^* \text{ is nonstationary in } \delta.$$

In particular, for $\zeta \geq \zeta_i$, the set $h_i^{-1}(\{\zeta\}) \cap S' \cap e_\delta$ is nonstationary in $\delta$.

Now let $\zeta^* \stackrel{\text{def}}{=} \sup\{\zeta_i : i < \theta\}$, hence $\zeta^* < \theta_1^+$ and

$$\zeta \geq \zeta^* \Longrightarrow (\forall i < \theta)[h_i^{-1}(\{\zeta\}) \cap S' \cap e_\delta \text{ is nonstationary in } \delta].$$

As $\text{cf}(\delta) = \theta_1^+$, there is $\alpha^* < \delta$ such that

$$\alpha \in e_\delta \setminus \alpha^* \Longrightarrow \text{otp}(e_\alpha) > \zeta^*.$$

Without loss of generality $\alpha^* > \theta_1^+$. For such $\alpha \in S'$, let

$$i(\alpha) \stackrel{\text{def}}{=} \min\{i < \theta : \text{otp}(e_{\xi_{\alpha,i}}) > \zeta^*\},$$

so $i(\alpha) < \alpha$. Now note that, as $\delta \in S_1$, we have that $S^*$ is stationary in $\delta$, so $e_\delta \cap S^*$ is stationary in $\delta$. On the other hand, $e_\delta \subseteq S_2$, hence $S' = S^* \cap S_2$ is stationary in $\delta$. Hence there is $T \subseteq S' \cap e_\delta$ stationary in $\delta$ and $i^* < \theta$ such that $\alpha \in T \Longrightarrow i(\alpha) = i^*$. As $\min(T) > \alpha^* > \theta_1^+$, the function $h_{i^*}$ is regressive on $T$, hence there is $T' \subseteq T$ stationary in $\delta$ and $\zeta^{**} > \zeta^*$ such that

$$\alpha \in T' \Longrightarrow h_{i^*}(\alpha) = \zeta^{**}.$$

Contradiction. $\bigstar_{1.3}$

Now let $i^*$ and $\beta^*$ be as guaranteed by Subclaim 1.3. Let $h \stackrel{\text{def}}{=} h_{i^*}^{\beta^*}$.

**Subclaim 1.4.** For some $\gamma = \gamma^* \leq \beta^*$ with $\text{cf}(\gamma) = \aleph_0$, for every club $E$ of $\lambda^+$, there are stationary many $\delta \in S_1 \cap S_2$ such that $e_\delta \subseteq E$ and

$$\gamma = \sup\{\zeta < \gamma : h^{-1}(\{\zeta\}) \cap S' \cap E \text{ is stationary in } \delta\}.$$

**Proof of the Subclaim.** Suppose not. Then for every $\gamma \leq \beta^*$ with $\mathrm{cf}(\gamma) = \aleph_0$, there are clubs $E_\gamma^0$ and $E_\gamma^1$ of $\lambda^+$ such that

$$[\delta \in S_1 \cap S_2 \cap E_\gamma^1 \ \& \ e_\delta \subseteq E_\gamma^0] \Longrightarrow$$
$$\gamma > \sup\{\zeta < \gamma : h^{-1}(\{\zeta\}) \cap S' \cap E_\gamma^0 \text{ is stationary in } \delta\}.$$

Let $E^* \stackrel{\mathrm{def}}{=} \mathrm{acc}(\bigcap_{\gamma \in S_{\aleph_0}^{\beta^*+1}} E_\gamma^0 \cap \bigcap_{\gamma \in S_{\aleph_0}^{\beta^*+1}} E_\gamma^1)$, so a club of $\lambda^+$. By the choice of $\beta^*$, there is $\delta \in E^* \cap S_1 \cap S_2$ such that $e_\delta \subseteq E^*$ and

$$|\{\zeta \leq \beta^* : h^{-1}(\{\zeta\}) \cap S' \cap E^* \text{ is stationary in } \delta\}| \geq \theta_1.$$

Hence we can choose $\{\gamma_n : n < \omega\}$ increasing, with $\gamma_n < \beta^*$ and such that $h^{-1}(\{\gamma_n\}) \cap S' \cap E^*$ is stationary in $\delta$, for all $n < \omega$. Let $\gamma \stackrel{\mathrm{def}}{=} \sup\{\gamma_n : n < \omega\}$, hence $\mathrm{cf}(\gamma) = \aleph_0$ and $\gamma \leq \beta^*$. This contradicts the choice of $E^*$. $\bigstar_{1.4}$

Now choose $\gamma^*$ as in Subclaim 1.4, and let $\langle \gamma_n^* : n < \omega \rangle$ be increasing to $\gamma^*$.

**Subclaim 1.5.** For some $a = a^* \in [\omega]^{\aleph_0}$, for every club $E$ of $\lambda^+$, there are stationarily many $\delta \in S_1 \cap S_2$ such that $e_\delta \subseteq E$ and

$$n \in a \Longrightarrow (\exists \zeta \in [\gamma_n^*, \gamma_{n+1}^*))[h^{-1}(\{\zeta\}) \cap S' \cap E^* \text{ stationary in } \delta].$$

**Proof of the Subclaim.** Suppose not. Then for every $a \in [\omega]^{\aleph_0}$, there are clubs $E_a^0$ and $E_a^1$ of $\lambda^+$, such that

$$[\delta \in S_1 \cap S_2 \cap E_a^1 \ \& \ e_\delta \subseteq E_a^0] \Longrightarrow$$
$$(\exists n_a \in a)(\forall \zeta \in [\gamma_{n_a}^*, \gamma_{n_a+1}^*))[h^{-1}(\{\zeta\}) \cap S' \cap E_a^0 \text{ is non-stationary in } \delta\}.$$

Let $E^* \stackrel{\mathrm{def}}{=} \mathrm{acc}(\bigcap_{a \in [\omega]^{\aleph_0}} E_a^0 \cap E_a^1)$, a club of $\lambda^+$, as $\lambda \geq 2^{\aleph_0}$. Hence for every $\delta \in S_1 \cap S_2 \cap E^*$ such that $e_\delta \subseteq E^*$ we have that for all $a \in [\omega]^{\aleph_0}$

$$(\exists n_a \in a)(\forall \zeta \in [\gamma_{n_a}^*, \gamma_{n_a+1}^*))[h^{-1}(\{\zeta\}) \cap S' \cap E^* \text{ is nonstationary in } \delta],$$

in other words

$$\{n : (\exists \zeta \in [\gamma_n^*, \gamma_{n+1}^*))[h^{-1}(\{\zeta\}) \cap S' \cap E^* \text{ is stationary in } \delta]\}$$

is bounded in $\omega$. Choosing a $\delta \in S_1 \cap S_2 \cap E^*$ such that $e_\delta \subseteq E^*$ and $\gamma^* = \sup\{\zeta < \gamma^* : h^{-1}(\{\zeta\}) \cap S' \cap E^*$ is stationary in $\delta\}$, we obtain a contradiction. $\bigstar_{1.5}$

Let $a^*$ be as in Subclaim 1.5, and let us enumerate $a^* \stackrel{\text{def}}{=} \{k_n : n < \omega\}$ increasingly. Let

$$g^{-1}(\{n\}) \stackrel{\text{def}}{=} \bigcup\{h^{-1}(\{\zeta\}) : \zeta \in [\gamma^*_{k_n}, \gamma^*_{k_{n+1}})\}.$$

$\bigstar_{1.2}$

We shall now define $\bar{c} = \langle c_\alpha : \alpha < \lambda^+\rangle$, so that

($\alpha$) For every $\alpha$, we have that $c_\alpha$ is a club of $\alpha$ with $\text{otp}(c_\alpha) \leq \lambda$,

($\beta$) If $\delta \in S_2$, then $c_\delta \supseteq e_\delta$,

($\gamma$) If $\delta \in S_2$ and $\sup(e_\delta) = \delta$, then $c_\delta = e_\delta$,

Now for any limit $\delta < \lambda^+$ we choose by induction on $n < \omega$ a club $C^n_\delta$ of $\delta$ of order type $\leq \lambda^{n+1}$, using the following algorithm:
Let $C^0_\delta \stackrel{\text{def}}{=} c_\delta$. Let

$$C^{n+1}_\delta \stackrel{\text{def}}{=} C^n_\delta \cup \{\alpha : (\exists \beta \in \text{nacc}(C^n_\delta)) [\sup(\beta \cap C^n_\delta) < \alpha < \beta \;\&\; \alpha \in c_\beta]\}.$$

**Note 1.6.** (1) The above algorithm really gives $C^n_\delta$ which is a club of $\delta$ with

$$\text{otp}(C^n_\delta) \leq \lambda^{n+1}.$$

If $\delta \in S_2$, then $\text{otp}(C^\delta_n) \leq \lambda^n \times \theta^+_1$.

[Why? We prove this by induction on $n$. It is clearly true for $n = 0$. Assume its truth for $n$. Clearly $C^{n+1}_\delta$ is unbounded in $\delta$, let us show that it is closed. Suppose $\alpha = \sup(C^{n+1}_\delta \cap \alpha) < \delta$. If $\alpha = \sup(C^n_\delta \cap \alpha)$, then $\alpha \in C^n_\delta \subseteq C^{n+1}_\delta$ by the induction hypothesis. So, assume

$$\alpha^* \stackrel{\text{def}}{=} \sup(C^n_\delta \cap \alpha) < \alpha$$

and $\alpha \notin C^n_\delta$. Let $\langle \alpha_i : i < \text{cf}(\alpha)\rangle$ be an increasing to $\alpha$ sequence in $(\alpha^*, \alpha) \cap C^{n+1}_\delta$. Hence for every $i$ there is $\beta_i \in \text{nacc}(C^n_\delta)$ such that $\alpha_i \in c_{\beta_i}$ and $\sup(C_\delta^n \cap \beta_i) < \alpha_i$. As $\sup(C_\delta^n \cap \alpha) = \alpha^* < \alpha_i$ and $\alpha \notin C_\delta^n$, we have $\beta_i > \alpha$, for every $i$. Suppose that $i \neq j$ and $\beta_i < \beta_j$. Hence $\sup(C_\delta^n \cap \beta_j) \geq \beta_i > \alpha_j$, a contradiction. So, there is $\beta$ such that $\beta_i = \beta$ for all $i$, hence $\{\alpha_i : i < \mathrm{cf}(\alpha)\} \subseteq c_\beta$. As $c_\beta$ is closed, we have $\alpha \in c_\beta$, and by the definition of $C_\delta^{n+1}$ we have $\alpha \in C_\delta^{n+1}$.

As for every $\beta$ we have $\mathrm{otp}(c_\beta) \leq \lambda$, and by the induction hypothesis $\mathrm{otp}(C_\delta^n) \leq \lambda^{n+1}$, we have $\mathrm{otp}(C_\delta^{n+1}) \leq \lambda^{n+2}$.

Similarly, if $\delta \in S_2$ clearly $\mathrm{otp}(C_\delta^n) \leq \lambda^n \times \theta_1^+$.]

(2) For every $n$, we have $\mathrm{acc}(C_\delta^n) \subseteq S_{<\lambda}^{\lambda^+}$.

[Why? Again by induction on $n$. For $n = 0$ it follows as $\mathrm{otp}(c_\delta) \leq \lambda$. Suppose this is true for $C_\delta^n$. The analysis from the proof of (1) shows that for $\alpha \in \mathrm{acc}(C_\delta^{n+1}) \setminus \mathrm{acc}(C_\delta^n)$, there is $\beta$ such that $\alpha \in c_\beta$, hence $\mathrm{cf}(\alpha) < \lambda$.]

(3) For every limit $\delta < \lambda^+$, we have $S_\lambda^\delta = \bigcup_{n<\omega} \mathrm{nacc}(C_\delta^n) \cap S_\lambda^\delta$.

[Why? Fix such $\delta$ and let $\alpha \in S_\lambda^\delta$. By item (2), it suffices to show that $\alpha \in C_\delta^n$ for some $n$. Suppose not, so let $\gamma_n \stackrel{\mathrm{def}}{=} \min(C_\delta^n \setminus \alpha)$ for $n < \omega$. Hence $\langle \gamma_n : n < \omega \rangle$ is a non-increasing sequence of ordinals $> \alpha$, and so there is $n^*$ such that $n \geq n^* \implies \gamma_n = \gamma_{n^*}$. In particular we have that $\gamma_{n^*} \in \mathrm{nacc}(C_\delta^{n^*})$. Let $\beta \in c_{\gamma_{n^*}} \setminus \alpha$. Hence $\sup(\beta \cap C_\delta^{n^*}) < \alpha \leq \beta < \gamma_{n^*}$. By the definition of $C_\delta^{n^*+1}$, we have $\beta \in C_\delta^{n^*+1}$, a contradiction.]

Now for each $\delta \in S_1 \cap S_2$ we define

$$E_\delta \stackrel{\mathrm{def}}{=} e_\delta \cup \bigcup \{C_\alpha^{g(\alpha)} \setminus \sup(e_\delta \cap \alpha) : \alpha \in \mathrm{nacc}(e_\delta) \cap S'\}.$$

Note first that $E_\delta$ is a club of $\delta$, for $\delta \in S'$.

[Why? Clearly, $E_\delta$ is unbounded. Suppose $\gamma = \sup(E_\delta \cap \gamma) < \delta$. Without loss of generality we can assume $\gamma \notin e_\delta$. Let $\gamma^* \stackrel{\mathrm{def}}{=} \sup(e_\delta \cap \gamma) < \gamma$. For every $\beta \in E_\delta \cap (\gamma^*, \gamma)$, there is $\alpha_\beta \in \mathrm{nacc}(e_\delta) \cap S'$ such that $\beta \in C_{\alpha_\beta}^{g(\alpha_\beta)} \setminus \sup(e_\delta \cap \alpha_\beta)$. By the choice of $\gamma^*$, every such $\alpha_\beta > \gamma$. Suppose that $\beta_1 \neq \beta_2 \in E_\delta \cap (\gamma^*, \gamma)$ and $\alpha_{\beta_1} < \alpha_{\beta_2}$. Hence $\sup(e_\delta \cap \alpha_{\beta_2}) \geq \alpha_{\beta_1}$, a contradiction. So all $\alpha_\beta$ are a fixed $\alpha$. Hence $\gamma < \alpha$ is a limit point of $C_\alpha^{g(\alpha)}$, and we are done, as $C_\alpha^{g(\alpha)}$ is closed.]

Also note that $\text{otp}(E_\delta) < \lambda^\omega \times \theta_1^+$, and that the above argument shows that $\text{acc}(E_\delta) \subseteq S_\lambda^{\lambda^+}$.

Suppose that $A \subseteq S_\lambda^{\lambda^+}$ is unbounded and it exemplifies that $\langle E_\delta : \delta \in S' \rangle$ fails to satisfy the requirements. Hence there is a club $E$ of $\lambda^+$ such that

$$\delta \in E \cap S' \implies \sup(A \cap \text{nacc}(E_\delta)) < \delta.$$

Let $E^* \stackrel{\text{def}}{=} \text{acc}(E) \cap \{\delta : \delta = \sup(A \cap \delta)\}$, hence a club of $\lambda^+$. Let $\delta^* \in S_1 \cap S_2 \cap E^*$ be such that $e_{\delta^*} \subseteq E^*$ and for all $n < \omega$, the set $E^* \cap \delta^* \cap g^{-1}(\{n\})$ is stationary in $\delta^*$.

For $\alpha \in \text{nacc}(e_{\delta^*})$ we have that $A \cap \alpha$ is unbounded in $\alpha$, hence by Note 1.6(3) (as $A \subseteq S_\lambda^{\lambda^+}$ and $\text{nacc}(e_{\delta^*}) \subseteq S_{\aleph_1}^{\lambda^+}$), there is $n < \omega$ such that $A \cap \text{nacc}(C_\alpha^n)$ is unbounded in $\alpha$. Let $n^*(\alpha)$ be the smallest such $n$. There is $n^*$ such that

$$\sup\{\alpha \in \text{nacc}(e_{\delta^*}) : n^*(\alpha) = n^*\} = \delta^*,$$

as $\text{cf}(\delta^*) > \aleph_0$. Let

$$e \stackrel{\text{def}}{=} \{\beta \in \text{acc}(e_{\delta^*}) : \beta = \sup\{\alpha \in \beta \cap \text{nacc}(e_{\delta^*}) : n^*(\alpha) = n^*\}\},$$

hence $e$ is a club of $\delta^*$. By our assumption, $\alpha^* \stackrel{\text{def}}{=} \sup(A \cap \text{nacc}(E_{\delta^*})) + 1 < \delta^*$. By the choice of $g$, the set $g^{-1}(\{n^*\}) \cap \delta^*$ is stationary in $\delta^*$. So, there is $\beta \in e \setminus \alpha^*$ such that $g(\beta) = n^*$. Let $\alpha \in (\alpha^*, \beta)$ be such that $\alpha \in \text{nacc}(e_{\delta^*})$ and $n^*(\alpha) = n^*$. Hence $A \cap \text{nacc}(C_\alpha^{n^*})$ is unbounded in $\alpha$. However,

$$C_\alpha^{n^*} \setminus \sup(\alpha \cap e_{\delta^*}) \subseteq E_{\delta^*},$$

hence there is $\gamma \in A \cap \text{nacc}(C_\beta^{n^*})$ with $\gamma \in E_{\delta^*} \setminus \alpha^*$. As $\gamma \in A$, we have $\text{cf}(\gamma) = \lambda$, hence $\gamma \in \text{nacc}(E_{\delta^*})$, a contradiction with the choice of $\alpha^*$.

(2) The assumption that $\lambda \geq 2^{\aleph_0}$ was used only in the proof of Claim 1.2. If $S^* = S_\theta^{\lambda^+}$, or just if $\{\text{otp}(e_\alpha) : \alpha \in e_\delta \cap S^* \ \& \ \alpha = \sup(e_\alpha)\}$ does not depend on $\delta$ of cofinality $\theta_1^+$ (which is true if $S^* = S_\theta^{\lambda^+}$), then the conclusion of Claim 1.2 easily follows. Namely, if $\delta \in S_1$, then $\text{cf}(\delta) = \theta_1^+ > \aleph_0$, and $\text{otp}(e_\delta) = \theta_1^+$. Hence, the set $C \stackrel{\text{def}}{=} \{\text{otp}(e_\alpha) : \alpha \in e_\delta \cap S^* \ \& \ \alpha = \sup(e_\alpha)\}$ is a club of $\theta_1^+$. Let $C = \bigcup_{n<\omega} T_n$, where each $T_n$ is a stationary subset of $\theta_1^+$ and $T_n$'s are pairwise disjoint. Define $g$ on $S'$ by letting $g(\alpha) = n$ if $\mathrm{otp}(e_\alpha) \in T_n$, and if $\alpha \in S'$ and $\mathrm{otp}(e_\alpha) \notin C$, just let $g(\alpha) = 0$. ★$_{1.1}$

## 2  A negation of guessing

**Theorem 2.1.** Assume that there is a supercompact cardinal. Then

(1) It is consistent that there is $\lambda$ a strong limit singular of cofinality $\omega$, such that $2^\lambda > \lambda^+$ and

(∗) There is a function $f : \lambda^+ \to \omega$ such that for every $\mathcal{P} \subseteq [\lambda^+]^\omega$ of cardinality $< 2^\lambda$, for some $X \in [\lambda^+]^{\lambda^+}$ we have

(i) $(\forall i < \omega)[|X \cap f^{-1}(\{i\})| = \lambda^+]$,
(ii) If $a \in \mathcal{P}$, then $\sup(\mathrm{Rang}(f \restriction (a \cap X))) < \omega$.

(2) In (1) we can replace $\omega$ by any regular $\kappa < \lambda$, but in the conclusion we do not necessarily obtain that $\lambda$ is a strong limit.

**Remark 2.2.** So the theorem basically states that no $\mathcal{P}$ as above provides a guessing.

**Proof.** (1) We start with a universe in which $\lambda$ is a supercompact cardinal and $GCH$ holds. We extend the universe by Laver's forcing ([La]), which makes the supercompactness of $\lambda$ indestructible by any extension by a $(<\lambda)$-directed-closed forcing. This forcing will preserve the fact that $2^\lambda = \lambda^+$. Let us call the so obtained universe $V$.

Now choose $\mu$ such that $\mu = \mu^\lambda > \lambda^+$. By [Ba], there is a $(<\lambda)$-closed $\lambda^{++}$-cc forcing notion $P$ adding $\mu$ unbounded subsets $A_\alpha$ ($\alpha < \mu$) to $\lambda^+$ such that

(∗∗) $\alpha \neq \beta < \mu \Longrightarrow |A_\alpha \cap A_\beta| < \lambda$.

In particular in $V^P$ we have $2^\lambda \le \mu$ and $\lambda$ is supercompact. In $V^P$, let $Q$ be Prikry's forcing which does not collapse cardinals and makes $\lambda$ singular with $\operatorname{cf}(\lambda) = \omega$, [Pr]. As this forcing does not add bounded subsets to $\lambda$, in the extension $\lambda$ is a strong limit singular and clearly satisfies $2^\lambda \le \mu$. In $V^{P*Q}$ we have (**). We now work in $V^{P*Q}$.

Let $\lambda = \sum_{\zeta < \omega} \lambda_\zeta$ where each $\lambda_\zeta < \lambda$ is regular. Let $\chi$ be large enough regular and $M \prec (\mathcal{H}(\chi), \in)$ with $\|M\| = \lambda^+$ be such that $\lambda^+ \subseteq M$ and $\langle A_\alpha : \alpha < \mu \rangle \in M$. We list $\bigcup_{\zeta < \omega}([\lambda^+]^{\lambda_\zeta} \cap M)$ as $\{b_i : i < \lambda^+\}$.

We define $f : \lambda^+ \to \omega$ by $f(i) = \zeta$ iff $|b_i| = \lambda_\zeta$. For $\alpha < \mu$, let $X_\alpha \stackrel{\text{def}}{=} \{i : b_i \subseteq A_\alpha\}$.

Now suppose that $\mathcal{P} \subseteq [\lambda^+]^\omega$ is of cardinality $< 2^\lambda \le \mu$, we shall look for $X$ as required in (*).

If $\alpha < \mu$ is such that $X_\alpha$ fails to serve as $X$, then one of the following two cases must hold:

<u>Case 1</u>. For some $\zeta < \omega$ we have $|\{i : b_i \subseteq A_\alpha \ \& \ |b_i| = \lambda_\zeta\}| < \lambda^+$, or

<u>Case 2</u>. For some $a \in \mathcal{P}$ we have $\sup(\operatorname{Rang}(f \restriction (a \cap X_\alpha))) = \omega$.

Considering the second case, we shall show that for any $a \in \mathcal{P}$, there are $< \lambda$ ordinals $\alpha$ such that the second case holds for $X_\alpha, a$. Fix an $a \in \mathcal{P}$. If $\alpha < \mu$ is such that Case 2 holds for $X_{\alpha,a}$, then

$$\sup(\{\zeta : (\exists i \in a)[b_i \subseteq A_\alpha \ \& \ |b_i| = \lambda_\zeta]\}) = \omega.$$

For $\zeta < \omega$ and $\alpha < \mu$ let $B_\zeta^\alpha \stackrel{\text{def}}{=} \{i \in a : b_i \subseteq A_\alpha \ \& \ |b_i| = \lambda_\zeta\}$. Notice that if $\alpha \ne \beta < \mu$ we have that for some $\zeta_{\alpha,\beta}$ the intersection $A_\alpha \cap A_\beta$ has size $< \lambda_{\zeta_{\alpha,\beta}}$, hence for all $\zeta \ge \zeta_{\alpha,\beta}$ we have $B_\zeta^\alpha \cap B_\zeta^\beta = \emptyset$.

Let $A \stackrel{\text{def}}{=} \{\alpha : \text{Case 2 holds for } a, \alpha\}$. For every $\alpha \in A$, let

$$S_\alpha \stackrel{\text{def}}{=} \{\zeta < \omega : B_\zeta^\alpha \ne \emptyset\},$$

hence $\alpha \ne \beta \implies S_\alpha \ne S_\beta$. Hence $|A| \le 2^\omega < \lambda$.

Now note that if $\alpha \in M$, then $A_\alpha \in M$, so for every $\zeta < \omega$ we have

$$\{i : b_i \subseteq A_\alpha \ \& \ |b_i| = \lambda_\zeta\} = [A_\alpha]^{\lambda_\zeta} \cap M.$$

As $M \models ``|[A_\alpha]^{\lambda_\zeta}| > \lambda$", we have $|\{i : b_i \subseteq A_\alpha \ \& \ |b_i| = \lambda_\zeta\}| = \lambda^+$. Hence Case 1 does not happen for this $\alpha$.

As we can find $\alpha \in M$ such that Case 2 does not happen, we finish.
(2) Use Magidor's forcing from [Ma] in place of Prikry's forcing. $\bigstar_{2.1}$